\documentclass[12pt,a4paper,twoside]{amsart}

\usepackage[ps2pdf,colorlinks=true,urlcolor=blue,
citecolor=red,linkcolor=blue,linktocpage,pdfpagelabels,bookmarksnumbered,bookmarksopen]{hyperref}

\newcommand{\eps}{\varepsilon}

\newcommand{\R}{\mathbb{R}}

\newcommand{\RN}{{\mathbb{R}^N}}
\newcommand{\RT}{{\mathbb{R}^3}}

\renewcommand{\b }{\beta }

\newcommand{\g }{\gamma }

\newcommand{\n }{\nabla }

\renewcommand{\P}{{\mathcal P}}

\newcommand{\SO}{\mathcal{S}}

\newcommand{\D }{{\mathcal D}^{1,2}(\RT)}

\renewcommand{\H}{H^{1}(\RT)}
\newcommand{\HR}{H^{1}_{r}(\RT)}
\newcommand{\HH}{\mathbb{H}}
\newcommand{\HHR}{\mathbb{H}_{r}}

\newcommand{\irt }{\int_{\RT}}

\def\bbm[#1]{\mbox{\boldmath $#1$}}
\newcommand{\I}{\mathrm{i}}

\newtheorem{theorem}{Theorem}[section]
\newtheorem{lemma}[theorem]{Lemma}

\newtheorem{proposition}[theorem]{Proposition}

\theoremstyle{definition}
\newtheorem{definition}[theorem]{Definition}
\newtheorem{remark}[theorem]{Remark}

\def\timestring{\begingroup
   \count0 = \time
   \divide\count0 by 60
   \count2 = \count0   
   \count4 = \time
   \multiply\count0 by 60
   \advance\count4 by -\count0   
   \ifnum\count4<10
      \toks1 = {0}%
   \else
      \toks1 = {}%
   \fi
   \ifnum\count2<12
      \toks0 = {a.m.}%
   \else
      \toks0 = {p.m.}%
      \advance\count2 by -12
   \fi
   \ifnum\count2=0
      \count2 = 12
   \fi
   \number\count2:\the\toks1 \number\count4 \thinspace \the\toks0
\endgroup}

\begin{document}

\title[Coupled nonlinear Schr\"odinger systems]{A note on coupled nonlinear Schr\"odinger systems under the effect of general nonlinearities}

\author{A. Pomponio}
\address{Dipartimento di Matematica, Politecnico di Bari, via Amendola 126/B, I-70125 Bari, Italy}
\email{a.pomponio@poliba.it}
\author{S. Secchi}
\address{Dipartimento di Matematica ed Applicazioni, Universit\`a di Milano--Bicocca, via R.~Cozzi 53, ed. U5, I--20125 Milano}
\email{simone.secchi@unimib.it}


\thanks{The second author was partially supported by the MIUR national project \textit{Metodi variazionali ed equazioni differenziali non lineari}, PRIN 2006.}

\date{\today\ \timestring}

\maketitle

\begin{abstract}
We prove the existence of radially symmetric ground--states for the system of Nonlinear Schr\"{o}dinger equations
\begin{equation*} 
\begin{cases}
-\Delta u+ u=f(u)+\b u v^2& \hbox{in }\RT,
\\
-\Delta v+ v=g(v)+\b u^2 v& \hbox{in }\RT,
\end{cases}
\end{equation*}
under very weak assumptions on the two nonlinearities $f$ and $g$. In particular, no ``Ambrosetti--Rabinowitz'' condition is required.
\end{abstract}

\section{Introduction}

In the last years, nonlinear Schr\"odinger systems have been widely investigated by several authors.
These systems are models for different physical phenomena: the propagation in birefringent
optical fibers, Kerr--like photorefractive media in optics, and Bose--Einstein
condensates. Roughly speaking, two ore more semilinear Schr\"{o}dinger equations like
\begin{equation} \label{eq:sing}
- \Delta u + a u=u^3 \qquad \text{\rm in $\RT$}
\end{equation}
are coupled together. Equation (\ref{eq:sing}) describes the propagation of pulse in a nonlinear optical fiber, and the existence of a unique (up to translation) least energy solution has been proved. It turns out that this \textit{ground state solution} is radially symmetric with respect to some point, positive and exponentially decaying together with its first derivatives at infinity.

Unluckily,  we know (see \cite{Ka}) that single-mode
optical fibers are not really ``single-mode'', but actually
bimodal due to the presence of birefringence. This birefringence
can deeply influence the way  an optical evolves during
the propagation along the fiber. Indeed, it can occur that the
linear birefringence makes a pulse split in two, while
nonlinear birefringent traps them together against splitting.
The evolution of two
orthogonal pulse envelopes in birefringent optical fibers is
governed (see \cite{M1,M2}) by the  nonlinear Schr\"odinger system
\begin{equation}\label{eq:M}
\begin{cases}
\I \displaystyle{\frac{\partial \phi}{\partial t}} +\displaystyle{\frac{\partial^2 \phi}{\partial x^2}}+|\phi|^2 \phi +\b |\psi|^2 \phi=0,
\\
\null \\
\I \displaystyle{\frac{\partial \psi}{\partial t}} +\displaystyle{\frac{\partial^2 \psi}{\partial x^2}}+|\psi|^2 \psi +\b |\phi|^2 \psi=0,
\end{cases}
\end{equation}
where $\b$ is a positive constant depending on the anisotropy of the
fibers. System (\ref{eq:M}) is also important for industrial applications in
fiber communications systems \cite{HK} and all-optical switching devices, see \cite{I}. If one looks for standing wave solutions of (\ref{eq:M}),
namely solutions of the form
\[
\phi(x,t)= e^{\I \omega_1^2 t} u(x) \quad {\rm and} \quad \psi(x,t)= e^{\I \omega_2^2 t} v(x),
\]
then \eqref{eq:M} becomes
\begin{equation}\label{eq:M'}
\begin{cases}
-\displaystyle{\frac{\partial^2 u}{\partial x^2}} + u =u^3 + \b v^2 u  & {\rm in} \; \R,
\\
\null \\
-\displaystyle{\frac{\partial^2 v}{\partial x^2}} + \omega^2 v =v^3 + \b u^2 v  & {\rm in} \; \R,
\end{cases}
\end{equation}
with $\omega^2={\omega_2^2}/{\omega_1^2}$. 
Other physical
phenomena, such as Kerr--like photorefractive media in optics,
are also described by \eqref{eq:M'}, see \cite{AA,CCMS}. As a word of caution, (\ref{eq:M'}) possesses the ``simple'' solutions of the form $(u,0)$ and $(0,v)$, where $u$ and $v$ solve (\ref{eq:sing}). 

Problem \eqref{eq:M'}, in a more general situation and also in higher dimension, has been studied in \cite{CZ1,CZ2}, where smooth ground state solutions $(u,v)\neq (0,0)$ are found
 by concentration compactness arguments. Later on, Ambrosetti \textit{et al.} in \cite{AC}, 
Maia \textit{et al.} in \cite{MMP} and Sirakov in \cite{S} deal with problem
\begin{equation}\label{eq:3}
\begin{cases}
-\Delta u+ u=u^3+\b u v^2& \hbox{in }\RT,
\\
-\Delta v+ v=v^3+\b u^2 v& \hbox{in }\RT,	
\end{cases}
\end{equation}
and, among other results, they prove the existence of
ground state solutions of the type $(u,v)$, with $u,v>0$, for $\b>0$ sufficiently big. Similar problems have been treated also in \cite{HS,LW1,LW3,WW}. Some results in the singularly perturbed case can be found in \cite{LW2,MPS,P}, while the orbital stability and blow-up proprieties have been studied in \cite{FM,MMP2}. 

Although the interest lies in solutions with both non-trivial components, solutions of \eqref{eq:3} are somehow related to solutions of the single nonlinear Schr\"odinger equation
(\ref{eq:sing}). The \emph{nonlinearity} $g(u)=u^3$ is typical in physical models, but much more general Schr\"odinger equation
of the form
\begin{equation*}
-\Delta u + a u=g(u), \quad \hbox{in }\RT,
\end{equation*}
still have at least a  ground state solution under general assumptions on the nonlinearity $g$ which, for example, do not require any Ambrosetti--Rabinowitz growth condition.  We recall that a function $f\colon \mathbb{R} \to [0,+\infty)$ satisfies the Ambrosetti--Rabinowitz growth condition if there exists some $\mu >2$ such that 
\begin{equation} \label{AR}
0< \mu F(t) \leq f(t)t \quad\hbox{for all $t \geq 0$},
\end{equation}
where $F(t)=\int_0^t f(s)\, ds$.
This condition essentially states that the function $\log \left( F(t)t^{-\mu} \right)$ is monotone increasing for $t \geq 0$, and is trivially satisfied by any power $f(u)=u^r$ with $r > 1$. The crucial feature is, to summarize, that ground states are necessarily radially symmetric with respect to some point, and this knowledge recovers some compactness. We refer to the celebrated papers \cite{BL1,CGM} for a deep study of these scalar--field equations (see also \cite{AP,JT}).

Motivated by these remarks, we want to find ground state solutions for the system
\begin{equation} \label{eq}
\begin{cases}
-\Delta u+ u=f(u)+\b u v^2& \hbox{in }\RT,
\\
-\Delta v+ v=g(v)+\b u^2 v& \hbox{in }\RT,
\end{cases}
\end{equation}
where $\b\in \R$ and $f,g\in C(\RT,\R)$ satisfy the following assumptions:
\begin{itemize}
	\item[{\bf (f1)}] $\displaystyle\lim_{t\to 0}\frac{f(t)}{t}=0$;
	\item[{\bf (f2)}] $\displaystyle\lim_{t\to \infty}\frac{f(t)}{|t|^p}=0$, for some $1<p<5$;
	\item[{\bf (f3)}] there exists $T_1>0$ such that $\frac 12 T_1^2< F(T_1)$, where $F(t)=\int_0^t f(s)\,d s$;
	\item[{\bf (g1)}] $\displaystyle\lim_{t\to 0}\frac{g(t)}{t}=0$;
	\item[{\bf (g2)}] $\displaystyle\lim_{t\to \infty}\frac{g(t)}{|t|^q}=0$, for some $1<q<5$;
	\item[{\bf (g3)}] there exists $T_2>0$ such that $\frac 12 T_2^2< G(T_2)$, where $G(t)=\int_0^t g(s)\,d s$.
\end{itemize}
\begin{remark}\label{re:AR}
These assumptions are very weak, and it is easy to construct function $f$ and $g$ that match them but do not match the Ambrosetti--Rabinowitz condition. Without condition (\ref{AR}), one cannot perform a standard minimization over the \emph{Nehari manifold} $\mathcal{N}=\{(u,v) \neq (0,0) \mid I'(u,v)[u,v]=0\}$ (see below for the definition of the functional $I$), as done in \cite{AC,S}. We recall that (\ref{AR}) is also used to prove the boundedness of minimizing sequences constrained to $\mathcal{N}$, so its failure can cause troubles even at this stage. 
Our existence results not only cover more general systems than those in \cite{AC,S}, but also give a different existence proof when $f$ and $g$ coincide with pure powers.
\end{remark}
\begin{remark}
For a single Schr\"{o}dinger equation, even weaker assumptions can be requested, see \cite{BL1}. Unluckily, the idea of locating solutions by means of the maximum principle does not seem to work for systems.
\end{remark}
System \eqref{eq} has a variational structure, in particular solutions of \eqref{eq} can be found as critical points of the functional $I:\H\times\H\to \R$ defined by
\begin{equation*}
I(u,v) = I_F(u)+I_G(v) - \frac{\b}{2} \int_{\RT} |u|^2 |v|^2\, dx,
\end{equation*}
where we have set 
\begin{align*}
I_F(u)&= \frac{1}{2} \|u\|_{H^1}^2 - \int_\RT F(u), \\
I_G(v) &= \frac{1}{2} \|v\|_{H^1}^2 - \int_\RT G(v).
\end{align*}

We will call \emph{ground state solution} any couple  $(u,v)\neq (0,0)$ which solves
(\ref{eq}) and minimizes the functional $I$ among all possible nontrivial solutions. 
Thus we have to overcome the strong lack of compactness under our weak assumptions on $f$ and $g$, and also to exclude ``simple'' solutions with a null component. To fix terminology, we introduce the following definition.

\begin{definition}
A solution of \eqref{eq}, $(u,v)\in \HH$, $(u,v)\neq(0,0)$ will be called {\it scalar solution} if either $u\equiv 0$ or $v\equiv 0$; while a solution of \eqref{eq}, $(u,v)\in \HH$, $(u,v)\neq(0,0)$ will be called {\it vector solution} if $u\neq0$ and $v\neq 0$.
\end{definition}

Scalar solutions for problem \eqref{eq} exist by the results of \cite{BL1}. Indeed, since $f$ satisfies ({\bf f1-3}),  there exists a (least--energy) solution $u_0 \in \H$ for the single Schr\"odinger equation
\begin{equation}	\label{eq:unaf}
-\Delta u+ u=f(u) \qquad \hbox{in }\RT,
\end{equation}
and since $g$ satisfies ({\bf g1-3}), there exists a (least--energy) solution $v_0 \in \H$ for
\begin{equation}	\label{eq:unag}
-\Delta v+ v=g(v) \qquad \hbox{in }\RT.
\end{equation}
It can be checked immediately that the couples $(u_0,0)$ and $(0,v_0)$ are non-trivial solutions of \eqref{eq}.

As a first step, we will prove that for any $\b\in \R$ the problem \eqref{eq} admits a ground state.

\begin{theorem}\label{main}
Let $f$ and $g$ satisfy ({\bf f1-3}) and ({\bf g1-3}). Then for any $\b\in \R$ there exists a ground state solution of \eqref{eq}. Moreover, if $\b>0$, this solution is radially symmetric.
\end{theorem}
Then we will prove that vector solutions exist whenever the coupling parameter $\beta$ is sufficiently large.
\begin{theorem}\label{main2}
Let $f$ and $g$ satisfy ({\bf f1-3}) and ({\bf g1-3}). Then there exists $\b_0>0$ such that, for any $\b>\b_0$, there there exists a vector solution of \eqref{eq}, which is a ground state solution. Moreover this solution is radially symmetric.
\end{theorem}

The main result of this paper is Theorem \ref{main2}: up to our knowledge, indeed, this is the first vector solution existence result for problem \eqref{eq}. As already said in Remark \ref{re:AR}, without condition (\ref{AR}), we cannot perform a standard minimization over the Nehari manifold and so our proof which is based on a constrained minimization over the Pohozaev manifold. The existence of vector solutions for any $\beta$ is not known. The existence result of Theorem \ref{main}, instead, is already known and it is proved in \cite{BL}. For the reader's sake, here we give a different proof based on the constrained minimization over the Pohozaev manifold.



\medskip

\begin{center}{\bf Notation}\end{center}
\begin{itemize} 
\item If $r>0$ and $x_0 \in \RT$, $B_r (x_0):= \left\{ x\in\RT : |x- x_0| <r \right\}$.
We denote with $B_r$ the ball of radius $r$ centered in the origin.
\item We denote with $\|\cdot \|$ the norm of $\H$.
\item We set $\HH=\H \times \H$ and, for any $(u,v)\in \HH$, we set $\|(u,v)\|^2=\|u\|^2+\|v\|^2$.
\item With $C_i$ and $c_i$, we denote generic positive constants, which may also vary from line to line.
\end{itemize}

\section{The Pohozaev manifold}

By ({\bf f1-2}) and ({\bf g1-2}), we get that for any $\eps>0$ there exists $C_\eps>0$ such that
\begin{align}
|f(t)|\leq \eps |t| +C_\eps |t|^p, & \qquad \hbox{for all }t \in \R;\label{eq:f}
\\
|F(t)|\leq \eps |t|^2 +C_\eps |t|^{p+1}, & \qquad \hbox{for all }t \in \R; \label{eq:F}
\\
|f(t)|\leq \eps |t| +C_\eps |t|^5, & \qquad \hbox{for all }t \in \R;\label{eq:f*}
\\
|F(t)|\leq \eps |t|^2 +C_\eps |t|^{6}, & \qquad \hbox{for all }t \in \R. \label{eq:F*} \\
|g(t)|\leq \eps |t| +C_\eps |t|^q, & \qquad \hbox{for all }t \in \R;\label{eq:g}
\\
|G(t)|\leq \eps |t|^2 +C_\eps |t|^{q+1}, & \qquad \hbox{for all }t \in \R; \label{eq:G}
\\
|g(t)|\leq \eps |t| +C_\eps |t|^5, & \qquad \hbox{for all }t \in \R;\label{eq:g*}
\\
|G(t)|\leq \eps |t|^2 +C_\eps |t|^{6}, & \qquad \hbox{for all }t \in \R. \label{eq:G*}
\end{align}

By \cite[Lemma 3.6]{ADM} and repeating the arguments of \cite{BL1}, it is easy to see that each solution of \eqref{eq}, $(u,v)\in \HH$, satisfies the following Pohozaev identity:
\begin{equation}	\label{eq:Poho}
\irt |\n u|^2+|\n v|^2
=6 \irt F(u)+G(v)-\frac{u^2}{2}  -\frac{v^2}{2}  +\frac \b 2 u^2 v^2.
\end{equation}
Therefore each non-trivial solution of \eqref{eq} belongs to $\P$, where
\begin{equation}	\label{eq:P}
\P:=\{(u,v)\in \HH \mid (u,v)\neq (0,0),\; (u,v) \hbox{ satisfies \eqref{eq:Poho} }\}.
\end{equation}
We call $\P$ the Pohozaev manifold associated to (\ref{eq}).
We collect its main properties of the set $\P$ in the next Proposition: the proof is easy and left to the reader.
\begin{proposition} \label{prop:P}
Define the functional $J \colon \mathbb{H} \to \mathbb{R}$ by
\begin{equation*}
J(u,v)=\frac 12\irt |\n u|^2+|\n v|^2 -3 \irt F(u)+G(v)-\frac{u^2}{2}  -\frac{v^2}{2}  +\frac \b 2 u^2 v^2.
\end{equation*}
Then 
\begin{enumerate}
\item $\P = \left\{ (u,v) \in \mathcal{H} \setminus \{(0,0)\} \mid J(u,v)=0 \right\}$;
\item $\P$ is a $C^1$--manifold of codimension one.
\end{enumerate}
\end{proposition}

\begin{lemma}\label{le:>C}
There exists $C>0$ such that $\|(u,v)\|\geq C$, for any $(u,v)\in \P$.
\end{lemma}

\begin{proof}
Let $(u,v)\in \P$. By \eqref{eq:F}, \eqref{eq:G} and \eqref{eq:P}, we easily get
\begin{align*}
\|u\|^2+\|v\|^2&\leq C_1 \irt |u|^{p+1} +|v|^{q+1} +u^2 v^2
\\
&\leq C_2 (\|u\|^{p+1}+\|v\|^{q+1}+ \|u\|^2\|v\|^2),
\end{align*}
which shows the claim.
\end{proof}

According to the definition of \cite{L1}, we say that a sequence $\{(u_n,v_n)\}_n$ vanishes if, for all $r>0$
\[
\lim_{n \to +\infty} \sup_{\xi \in \RT} \int_{B_r(\xi)}u_n^2+v_n^2=0.
\]
\begin{lemma}\label{le:nonvan}
Any bounded sequence $\{(u_n,v_n)\}_n\subset \P$ does not vanish.
\end{lemma}

\begin{proof}
Suppose by contradiction that $\{(u_n,v_n)\}_n$ vanishes, then, in particular there exists $\bar r>0$ such that
\[
\lim_{n \to +\infty} \sup_{\xi \in \RT} \int_{B_{\bar r}(\xi)}u_n^2=0,
\quad \lim_{n \to +\infty} \sup_{\xi \in \RT} \int_{B_{\bar r}(\xi)}v_n^2=0.
\]
Then, by \cite[Lemma 1.1]{L2}, we infer that $u_n,v_n\to 0$ in $L^s(\RT)$, for any $2<s<6$. Since $\{(u_n,v_n)\}_n\subset \P$, we have that $(u_n,v_n) \to 0$ in $\HH$, contradicting Lemma \ref{le:>C}.
\end{proof}


\begin{lemma}\label{le:nat}
For any $\b\in \R$, $\P$ is a natural constraint for the functional $I$.
\end{lemma}

\begin{proof}
First we show that the manifold is nondegenerate in the
following sense:
\begin{equation*}
J'(u,v)\neq 0 \quad \hbox{for all $(u,v)\in\P$}.
\end{equation*}
By contradiction, suppose that $(u,v)\in\P$ and $J'(u,v)= 0,$
namely $(u,v)$ is a solution of the equation
\begin{equation}\label{eq:po}
\left\{
\begin{array}{ll}
-\Delta u+ 3u=3f(u)+3\b u v^2& \hbox{in }\RT,
\\
-\Delta v+ 3v=3g(v)+3\b u^2 v& \hbox{in }\RT.	
\end{array}
\right.
\end{equation}
As a consequence, $(u,v)$ satisfies the Pohozaev identity referred to
\eqref{eq:po}, that is
\begin{equation}\label{eq:pohopoho}
\irt |\n u|^2+|\n v|^2
=18 \irt F(u)+G(v)-\frac{u^2}2 -\frac{v^2}2 +\frac \b 2 u^2 v^2.
\end{equation}
Since $(u,v)\in\P$, by \eqref{eq:pohopoho} we get
\begin{equation*}
2\irt |\n u|^2+|\n v|^2=0
\end{equation*}
and we conclude that $u=v=0$: we get a contradiction since $(u,v)\in\P.$
\\
Now we pass to prove that $\P$ is a natural constraint for $I$.
Suppose that $(u,v)\in\P$ is a critical point of the
functional $I_{|\P}.$ Then, by Proposition \ref{prop:P}, there exists $\mu\in\R$ such
that
\begin{equation*}
I'(u,v) = \mu J'(u,v).
\end{equation*}
As a consequence, $(u,v)$ satisfies the following Pohozaev
identity
\begin{equation*}
\mu^{-1} \,  J(u,v) = \frac 12\irt |\n u|^2+|\n v|^2 -9 \irt F(u)+G(v)-\frac{u^2}2 -\frac{v^2}2 +\frac \b 2 u^2 v^2
\end{equation*}
which, since $J(u,v)=0$, can be written
\begin{equation*}
\mu \irt |\n u|^2+|\n v|^2 = 0.
\end{equation*}
Since either $u\neq 0$ or $v\neq 0$ we deduce that $\mu=0$, and we conclude.
\end{proof}
We set
\[
m=\inf_{(u,v)\in \P}I(u,v).
\]
We set $\HHR=H^{1}_{r}(\RT)\times H^{1}_{r}(\RT)$: here $H^{1}_{r}(\RT)$ denotes the radially symmetric functions of $\H$.

By means of the previous lemma, we are reduced to look for a
minimizer of $I$ restricted to $\P.$ By the well known properties
of the Schwarz symmetrization, we are allowed to work on the
functional space $\HHR$ as shown by the following

\begin{lemma}\label{le:simm}
For any $\b>0$ and for any $(u,v)\in\P$, there exists $(\bar u, \bar v)\in\P\cap \HHR$ such
that $I(\bar u, \bar v)\leq I(u,v)$.
\end{lemma}
\begin{proof}
Let $(u,v)\in\P$ and set $u^*, v^*\in H^{1}_{r}(\RT)$ their respective symmetrized functions.
We have
\begin{align*}
\irt |\n u^*|^2+|\n v^*|^2
&\leq \irt |\n u|^2+|\n v|^2
\\
&=6 \irt F(u)+G(v)-\frac{u^2}{2}  -\frac{v^2}{2}  +\frac \b 2 u^2 v^2
\\
&\leq 6 \irt F(u^*)+G(v^*)-\frac{(u^*)^2}{2}  -\frac{(v^*)^2}{2}  +\frac \b 2 (u^*)^2 (v^*)^2.
\end{align*}
Hence, there exists $\bar t\in (0,1]$ such that
$(\bar u, \bar v):=(u^*(\cdot/\bar t),v^*(\cdot/\bar t))\in\P\cap\HHR$ and
\begin{align*}
I(\bar u,\bar v)&=\frac 1 3 \irt | \n \bar u|^2 +|\n \bar v|^2
= \frac {\bar t} 3 \irt | \n u^*|^2 +|\n v^*|^2
\\
&\leq \frac 1 3 \irt |\n u|^2 +|\n v|^2= I(u,v).
\end{align*}
\end{proof}

\begin{proposition}\label{pr:m}
For any $\b>0$, the value $m$ is achieved as a minimum by $I$ on $\P$ by $(u,v)\in \HHR$.
\end{proposition}

\begin{proof}
For any $(u,v)\in \P$ we have
\begin{equation}	\label{eq:IP}
I(u,v)=\frac 13 \irt |\n u|^2+|\n v|^2\geq 0.
\end{equation}
Let $\{(u_n,v_n)\}_n\subset \P$ be such that $I(u_n,v_n)\to m$. By Lemma \ref{le:simm}, we can assume that $\{(u_n,v_n)\}_n\subset \P\cap\HHR$. 
\\
By \eqref{eq:IP}, we infer that $\{u_n\}_n, \{v_n\}_n$ are bounded in $\D$.
\\
Let $\varepsilon>0$ be given, and let $C_\varepsilon>0$ be the positive constant as in (\ref{eq:f})--(\ref{eq:G*}).
We observe that, for any $\b\in \R$, there exists a positive constant $C>0$ such that for all $x,y\in \R$:
\[
\eps(x^2+y^2)+C_\eps(x^6+y^6)+\frac \b 2 x^2 y^2
\leq 2\eps(x^2+y^2)+C C_\eps(x^6+y^6).
\]
Hence, since $\{(u_n,v_n)\}_n\subset \P$, by \eqref{eq:F*} and \eqref{eq:G*}, we get
\begin{align*}
\|u_n\|^2 + \|v_n\|^2
&\leq 6\irt \eps(u_n^2 +v_n^2) +C_\eps (u_n^6 +v_n^6)+\frac \b 2 u_n^2 v_n^2
\\
&\leq 6\irt 2\eps(u_n^2 +v_n^2) +C C_\eps (u_n^6 +v_n^6).
\end{align*}
Since $\{u_n\}_n, \{v_n\}_n$ are bounded in $\D$ and $\D\hookrightarrow L^6(\RT)$, then $\{(u_n,v_n)\}_n$ is bounded in $\HH$.
\\
By Lemma \ref{le:nonvan} we know that $\{(u_n,v_n)\}_n$ does not vanish, namely there exist $C,r>0$, $\{\xi_n\}_n\subset \RT$ such that
\begin{equation}\label{eq:nonvan2}
\int_{B_r(\xi_n)} u_n^2 +v_n^2\geq C, \hbox{ for all }n\geq 1.
\end{equation}
Since we are dealing with radially symmetric functions, without loss of generality, we can assume that $\xi_n=0$, for all $n \geq 1$.
\\
Since $\{(u_n,v_n)\}_n$ is bounded in $\HHR$, there exist $u,v\in \HR$ such that, up to a subsequence,
\begin{align*}
u_n \rightharpoonup u \; \hbox{in }\HR; & \;\;\;\;
u_n \to u \; \hbox{a.e. in }\RT; &\!
u_n \to u \; \hbox{in } L^s(\RT), \;2<s<6;
\\
v_n \rightharpoonup v \; \hbox{in }\HR; &\;\;\;\;
v_n \to v \; \hbox{a.e. in }\RT; &\!
v_n \to v \; \hbox{in } L^s(\RT), \;2< s<6.
\end{align*}
By \eqref{eq:nonvan2}, we can argue that either $u\neq 0$ or $v\neq 0$ and, moreover, since $\{(u_n,v_n)\}_n \subset \P$, passing to the limit, we have
\begin{equation}	\label{eq:<Poho}
\irt \frac 12|\n u|^2+\frac 12|\n v|^2+\frac 32 u^2 +\frac 32 v^2
\leq 3 \irt F(u)+G(v) +\frac \b 2 u^2 v^2.
\end{equation}
By \eqref{eq:<Poho}, it is easy to see that there exists $\bar t\in (0,1]$ such that $(\bar u,\bar v)=(u(\cdot/\bar t),v(\cdot/\bar t))\in \P\cap \HHR$. By the weak lower semicontinuity, we get
\begin{align*}
b &\leq I(\bar u,\bar v)
=\frac{\bar t}{3} \irt I(u,v)
\leq \frac{1}{3} \irt |\n u|^2+|\n v|^2
\\
& \leq \liminf_{n \to +\infty} \frac 13 \irt |\n u_n|^2+|\n v_n|^2
= \liminf_{n \to +\infty} I(u_n,v_n)=b,
\end{align*}
hence $(\bar u,\bar v)$ is a minimum of $I$ restricted on $\P$ and so, by Lemma \ref{le:nat}, it is a (radially symmetric) ground state solution for the problem \eqref{eq}.
\end{proof}

Finally, let us prove a lemma which will be a key point in the proof of Theorem \ref{main2}

\begin{lemma}\label{le:gamma}
Let $u_0,v_0\in \H$ be two non-trivial solutions respectively of \eqref{eq:unaf} and \eqref{eq:unag}. Then, for any $\b>0$, there exists $\bar t>0$ such that
$(u_0(\cdot/\bar t),v_0(\cdot/\bar t)) \in \P$.
\end{lemma}

\begin{proof}
Since $u_0$ is a solution of \eqref{eq:unaf}, then it satisfies the following Pohozaev identity:
\[
\irt |\n u_0|^2+3 \irt u_0^2= 6\irt F(u_0),
\]
hence
\begin{equation}	\label{eq:<F}
\irt u_0^2<2\irt F(u_0).
\end{equation}
Analogously, $v_0$ satisfies
\begin{equation}	\label{eq:<G}
\irt v_0^2<2\irt G(v_0).
\end{equation}
We set $\g(t)=(u_0(\cdot/t),v_0(\cdot/t))$, with $t>0$. We have
\[
I(\g(t))=
\frac t2 \irt |\n u_0|^2 + |\n v_0|^2
+ t^3 \!\!\irt \frac{u_0^2}{2} + \frac{v_0^2}{2}  - F(u_0) -G(v_0)-\frac \b 2 u_0^2 v_0^2.
\]
Since $I(\g(t))>0$ for small $t$ and, by \eqref{eq:<F} and \eqref{eq:<G} and being $\b>0$, $\lim_{t\to +\infty}I(\g(t))=-\infty$, there exists $\bar t>0$ such that $\frac{d }{d t}I(\g(\bar t)) =0$, which implies that the couple $(u_0(\cdot/\bar t),v_0(\cdot/\bar t)) \in \P$.
\end{proof}

\section{Proofs of Theorem \ref{main} and \ref{main2}}

\textit{Proof of Theorem \ref{main}}. 
If we set
\begin{equation*}
\SO:=\left\{(u,v)\in \HH\mid (u,v)\neq (0,0), \; (u,v)\hbox{ solves~ \eqref{eq}}\right\},
\end{equation*}
\begin{equation*}
b:=\inf_{(u,v)\in \SO} I(u,v)\geq 0.
\end{equation*}

Let $\{(u_n,v_n)\}_n\subset \SO$ be such that $I(u_n,v_n)\to b$. By \eqref{eq:IP}, we infer that $\{u_n\}_n, \{v_n\}_n$ are bounded in $\D$.
\\
Repeating the arguments of the proof of Proposition \ref{pr:m}, we can argue that $\{(u_n,v_n)\}_n$ is bounded in $\HHR$.
\\
By Lemma \ref{le:nonvan} we know that $\{(u_n,v_n)\}_n$ does not vanish, namely there exist $C,r>0$, $\{\xi_n\}_n\subset \RT$ such that
\begin{equation}\label{eq:nonvan}
\int_{B_r(\xi_n)} u_n^2 +v_n^2\geq C, \hbox{ for all }n\geq 1.
\end{equation}
Due to the invariance by translations, without loss of generality, we can assume that $\xi_n=0$ for every $n$.

Since $\{(u_n,v_n)\}_n$ is bounded in $\HH$, there exist $u,v\in \H$ such that, up to a subsequence,
\begin{align*}
u_n \rightharpoonup u \; \hbox{in }\H; & \;\;\;\;
u_n \to u \; \hbox{a.e. in }\RT; &\!
u_n \to u \; \hbox{in } L^s_{\mathrm{loc}}(\RT), \;1\leq s<6;
\\
v_n \rightharpoonup v \; \hbox{in }\H; &\;\;\;\;
v_n \to v \; \hbox{a.e. in }\RT; &\!
v_n \to v \; \hbox{in } L^s_{\mathrm{loc}}(\RT), \;1\leq s<6.
\end{align*}
By \eqref{eq:nonvan}, we can argue that either $u\neq 0$ or $v\neq 0$ and then it is easy to see that $(u,v)\in \SO$. By the weak lower semicontinuity, we get
\begin{align*}
b &\leq I(u,v)
=\frac 13 \irt |\n u|^2+|\n v|^2
\\
& \leq \liminf_{n \to +\infty} \frac 13 \irt |\n u_n|^2+|\n v_n|^2
= \liminf_{n \to +\infty} I(u_n,v_n)=b,
\end{align*}
hence $(u,v)$ is a ground state for the problem \eqref{eq}.

If $\b>0$, by  Proposition \ref{pr:m}, we can argue that there exists a ground state $(u,v)$ which belongs to $\HHR$.
$\hfill\square$



\vspace{1cm}

\textit{Proof of Theorem \ref{main2}}. We will use some ideas from \cite{MMP}.
Let $u_0,v_0\in \HR$ be two ground state solutions respectively for equation \eqref{eq:unaf} and equation \eqref{eq:unag}. By Lemma \ref{le:gamma}, we know that there exists $\bar t>0$ such that $(u_0(\cdot/\bar t),v_0(\cdot/\bar t)) \in \P$. Then, to show that any radial ground state solution $(\bar u,\bar v)$ is a vector solution, it is sufficient to prove that, for $\b$ positive and sufficiently large,
\begin{equation}	\label{eq:II}
I(u_0(\cdot/\bar t),v_0(\cdot/\bar t))< \min\{I(u_0,0),I(0,v_0)\}.
\end{equation}
Indeed, with some calculations, we have
\[
I(u_0(\cdot/\bar t),v_0(\cdot/\bar t))=
\frac{\left(\frac 13 \irt  |\n u_0|^2+|\n v_0|^2\right)^{3/2}}
{\left(2\irt F(u_0)+G(v_0)+\frac \b 2 u_0^2 v_0^2 - \frac{u_0^2}{2} -\frac{v_0^2}{2} \right)^{1/2}}.
\]
Then, for $\b$ positive and sufficiently large, we have \eqref{eq:II}.
$\hfill\square$

\medskip

\noindent \textbf{Note added in proof.} After the manuscript had been submitted, the preprint \cite{BJM} was pointed out to us. It is shown, in particular, that the radial symmetry of solutions to our system can be proved \textit{a priori}. More precisely, \emph{every} ground-state solution is necessarily invariant under rotations.


\begin{thebibliography}{99}



\bibitem{AA}
N. Akhmediev, A. Ankiewicz, 
{\it Partially coherent solitons on a finite background}, 
Phys. Rev. Lett., {\bf 82}, (1999), 2661--2664.


\bibitem{AC}
A. Ambrosetti, E. Colorado, {\it Standing waves of some coupled nonlinear Schr\"odinger equations}, 
J. Lond. Math. Soc., {\bf 75}, (2007),  67--82.


\bibitem{ADM}
P. Amster, P. De N\'apoli, M.C. Mariani, {\it Existence of solutions for elliptic systems with critical Sobolev exponent}, Electron. J. Differential Equations, {\bf 2002}, No. 49, 13 pp. (electronic).


\bibitem{AP}
A. Azzollini, A. Pomponio {\it On the Schr\"odinger equation in
$\RN$ under the effect of a general nonlinear term}, Indiana Univ.
Journal, to appear.


\bibitem{BL1}
H. Berestycki, P.L. Lions, {\it Nonlinear scalar field equations,
I - Existence of a ground state}, Arch. Rational Mech. Anal., {\bf
82}, (1983), 313--345.



\bibitem{BL}
H. Brezis, E. Lieb, \textit{Minimum action solutions of some vector field equations}, Comm. Math. Phys. \textbf{96} (1984), 97--113. 

\bibitem{BJM}
J. Byeon, L. Jeanjean, M. Maris, \textit{Symmetry and monotonicity of least energy solutions}, Calc. Var. Partial Diff. Equations, in press.

\bibitem{CCMS}
D.N. Christodoulides, T.H. Coskun, M. Mitchell, M. Segev, 
{\it Theory of incoherent self-focusing in biased photorefractive media}, 
Phys. Rev. Lett., {\bf 78}, (1997), 646-649.

 
\bibitem{CZ1}
R. Cipolatti, W. Zumpichiatti, 
{\it On the existence and regularity of ground states for a nonlinear system of coupled 
Schr\"odinger equations}, 
Comput. Appl. Math., {\bf 18}, (1999), 19--36.


\bibitem{CZ2}
R. Cipolatti, W. Zumpichiatti,
{\it Orbitally stable standing waves for a system of coupled nonlinear Schr\"odinger equations}, 
Nonlinear Anal., {\bf 42}, (2000), 445--461.


\bibitem {CGM}
S. Coleman, V. Glaser, A. Martin, {\it Action minima amoung
solutions to a class of euclidean scalar field equations}, Commun.
math. Phys., {\bf 58}, (1978), 211--221.
%
%


\bibitem{FM}
L. Fanelli, E. Montefusco, {\it On the blow-up threshold for weakly coupled nonlinear Schr\"odinger equations}, 
J. Phys. A: Math. Theor., {\bf 40}, (2007), 14139--14150.




\bibitem{HK}
A. Hasegawa, Y. Kodama,
Solitions in optical communications, Academic Press, San Diego, 1995.


\bibitem{HS}
N. Hirano, N. Shioji, {\it Multiple existence of solutions for coupled nonlinear Schr\"odinger equations}, 
Nonlinear Anal., {\bf 68}, (2008), 3845--3859.


\bibitem{I}
M.N. Islam, 
Ultrafast fiber switching devices and systems, Cambridge University Press, New 
York, 1992.



\bibitem{JT}
L. Jeanjean, K. Tanaka, {\it A positive solution for a nonlinear
Schr\"odinger equation on $\RN$}, Indiana Univ. Math. J., {\bf
54}, (2005), 443--464.


\bibitem{Ka}
I.P. Kaminow,
{\it Polarization in optical fibers},
IEEE J. Quantum Electron., {\bf 17}, (1981), 15--22.




\bibitem{LW1}
T.C. Lin, J. Wei,
{\it  Ground State of $N$ Coupled Nonlinear Schr\"odinger Equations in $\R^n$, $n\geq 3$},
Comm. Math. Phys., {\bf  255}, (2005), 629--653.


\bibitem{LW2}
T.C. Lin, J. Wei,
{\it  Spikes in two coupled nonlinear Schr\"odinger equations},
Ann. Inst. H. Poincar\'e Anal. Non lin\'earie, {\bf 22}, (2005), 403--439.


\bibitem{LW3}
T.C. Lin, J. Wei, {\it Symbiotic bright solitary wave solutions of coupled nonlinear Schr\"odinger equations},  Nonlinearity, {\bf 19}, (2006), 2755--2773. 

                                         
\bibitem{L1}
P.L. Lions,
\textit{The concentration-compactness principle in the calculus of variation.
The locally compact case. Part I},
Ann. Inst. Henri Poincar\'e, Anal. Non Lin\'eaire, {\bf 1}, (1984), 109--145.


\bibitem{L2}
P.L. Lions,
\textit{The concentration-compactness principle in the calculus of variation.
The locally compact case. Part II},
Ann. Inst. Henri Poincar\'e, Anal. Non Lin\'eaire, {\bf 1}, (1984), 223--283.


\bibitem{MMP}
L.A. Maia, E. Montefusco, B. Pellacci,
{\it Positive solutions for a weakly coupled nonlinear Schr\"odinger system},
J. Differential Equations, {\bf 229}, (2006), 743--767.


\bibitem{MMP2}
L.A. Maia, E. Montefusco, B. Pellacci, {\it Orbital stability of ground state solutions of coupled nonlinear Schr\"odinger equations}, preprint.


\bibitem{M1}
C.R. Menyuk,
{\it Nonlinear pulse propagation in birefringence optical fiber},
IEEE J. Quantum Electron., {\bf 23}, (1987), 174--176.


\bibitem{M2}
C.R. Menyuk,
{\it Pulse propagation in an elliptically birefringent Kerr medium},
IEEE J. Quantum Electron., {\bf 25}, (1989), 2674--2682.


\bibitem{MPS}
E. Montefusco, B. Pellaccci, M. Squassina, {\it Semiclassical states for weakly coupled nonlinear Schr\"odinger syste}, 
J. Eur. Math. Soc., {\bf 10}, (2008), 47--71.


\bibitem{P}
A. Pomponio, {\it Coupled nonlinear Schr\"odinger systems with potentials}, 
J. Differential Equations, {\bf 227}, (2006), 258--281.


\bibitem{S}
B. Sirakov, {\it Least energy solitary waves for a system of nonlinear Schr\"odinger equations in $\R^n$}, 
Comm. Math. Phys., {\bf 271}, (2007), 199--221.


\bibitem{WW}
J. Wei, T. Weth, {\it Radial solutions and phase separation in a system of two coupled Schr\"odinger equations}, 
Arch. Ration. Mech. Anal., {\bf 190}, (2008), 83--106. 

\end{thebibliography}
\end{document}